# Calculating the Fundamental Group of the Circle in Homotopy Type Theory

Daniel R. Licata
Institute for Advanced Study
drl@cs.cmu.edu

Michael Shulman
Institute for Advanced Study
mshulman@ias.edu

*Abstract*—Recent work on homotopy type theory exploits an exciting new correspondence between Martin-Löf's dependent type theory and the mathematical disciplines of category theory and homotopy theory. The category theory and homotopy theory suggest new principles to add to type theory, and type theory can be used in novel ways to formalize these areas of mathematics. In this paper, we formalize a basic result in algebraic topology, that the fundamental group of the circle is the integers. Though simple, this example is interesting for several reasons: it illustrates the new principles in homotopy type theory; it mixes ideas from traditional homotopy-theoretic proofs of the result with type-theoretic inductive reasoning; and it provides a context for understanding an existing puzzle in type theory—that a universe (type of types) is necessary to prove that the constructors of inductive types are disjoint and injective.

## I. INTRODUCTION

Recently, researchers have discovered an exciting new correspondence between Martin-Löf's dependent type theory and the mathematical disciplines of category theory and homotopy theory [1, 3, 4, 6, 10, 11, 19, 20, 21]. Under this correspondence, a type $A$ in dependent type theory carries the structure of an ∞-groupoid, or a topological space up to homotopy. Terms $M : A$ correspond to objects of the groupoid, or points in the topological space. Terms of Martin-Löf's intensional identity type, written $\alpha : \mathsf{Id}_A(M,N)$, correspond to morphisms, or paths in the topological space, between $M$ and $N$. Iterating the identity type gives further structure; for example, the type $\mathsf{Id}_{\mathsf{Id}_A(M,N)}(\alpha,\beta)$ represents higher-dimensional morphisms, or homotopies between paths. This correspondence has many applications: The category theory and homotopy theory suggest new principles to add to type theory, such as higher-dimensional inductive types [12, 13, 17] and Voevodsky's univalence axiom [7, 20]. Proof assistants such as Coq [2] and Agda [14], especially when extended with these new principles, can be used in novel ways to formalize category theory, homotopy theory, and mathematics in general.

Here, we consider the use of type theory for computer-checked proofs in homotopy theory. Rather than working with some concrete implementation of homotopy types (such as topological spaces or simplicial sets), we use type theory to give an abstract, combinatorial description of them. In this way, type theory serves as a *logic of homotopy theory*. To illustrate this, we compute what is called the *fundamental group* of the circle. To explain the meaning of this, consider a topological space $X$. Given a particular *base point* $x_0 \in X$, the *loops* in $X$ are the continuous paths from $x_0$ to itself. These loops (considered up to homotopy) have the structure of a group: there is an identity path (standing still), composition (go along one loop and then another), and inverses (go backwards along a loop). Thus the homotopy-equivalence classes of such loops form a group called the fundamental group of $X$ at $x_0$, denoted $\pi_1(X, x_0)$ or just $\pi_1(X)$. More generally, by considering higher-dimensional paths and deformations, one obtains the *higher homotopy groups* $\pi_n(X)$. Characterizing these is a central question in homotopy theory; they are surprisingly complex even for a space as simple as the sphere.

Consider the circle (written $\mathsf{S}^1$) with some fixed base point base. What paths are there from base to base? One possibility is to stand still. Others are to go around clockwise once, or to go around clockwise twice, etc. Or we can go around counter-clockwise once, twice, etc. However, up to homotopy, going around clockwise and then counterclockwise (or vice versa) is the identity: we can deform this path continuously back to the constant one. Thus, the clockwise and counterclockwise paths are inverses in $\pi_1(\mathsf{S}^1)$. This suggests that $\pi_1(\mathsf{S}^1)$ should be isomorphic to $\mathbb{Z}$, the additive group of the on the integers: one can stand still (0), or go around counterclockwise $n$ times ($+n$), or go around clockwise $n$ times ($-n$). Proving this formally is one of the first basic theorems of algebraic topology.

In this paper, we formalize such a proof in type theory, using Agda. Though simple, this example is interesting for several reasons. First, it illustrates the new ingredients in homotopy type theory: spaces-up-to-homotopy can be described in a direct, logical way, which captures their (higher-dimensional) inductive nature. In particular, our "circle" has a direct inductive presentation rather than a topological one. Voevodsky's univalence axiom also plays an essential role in the proof. Second, as we discuss below, the development of this proof was an interplay between homotopy theory and type theory, mixing ideas from traditional homotopy-theoretic proofs with techniques that are common in type theory. Third, the proof has computational content: it can also be seen as a program that converts a path on the circle to its winding number, and vice versa. Finally, it provides a context for understanding the familiar puzzle that a universe (type of types) is necessary to prove seemingly obvious properties of inductive types, such as injectivity and disjointness of constructors.

The remainder of this paper is organized as follows. In Section II, we introduce the basic definitions of homotopy type theory. In Section III, we introduce the methodology that we will use to calculate $\pi_1(S^1)$ using a warm-up example, proving injectivity and disjointness for the constructors of the coproduct type. In Section IV, we define the circle as a higher-dimensional inductive type, and in Section V we prove that its fundamental group is $\mathbb{Z}$.

## II. BASICS OF HOMOTOPY TYPE THEORY

### A. Types as Spaces

The central observation in homotopy type theory is that Martin-Löf's intensional identity type (traditionally called *propositional equality*) equips each type with the structure of a *homotopy type*, or, equivalently, an *infinite-dimensional groupoid*. Homotopy types are an abstract structure, with many concrete realizations such as topological spaces up to homotopy and simplicial sets. To emphasize the homotopy-theoretic interpretation, we write the identity type between elements $M, N : A$ as Path $M\ N$ (leaving $A$ as an implicit argument[1]). Path is defined as an inductive family of types, with one constructor id:

```
data Path {A : Type} : A → A → Type where
    id : {M : A} → Path M M
```

id (reflexivity of propositional equality) represents the identity path in the type A. We use the identifier Type for what is normally called Set in Agda (the type of smaller types), because the word "set" has another meaning in this context.[2]

The standard *J* elimination rule for the identity type (in the Paulin-Mohring "one-sided" form [16]) expresses that the paths from a fixed point M, to a variable endpoint x, are inductively generated by M and id. We call this *path induction*:

```
path-induction : {A : Type} {M : A}
    (C : (x : A) → Path M x → Type) (b : C M id)
    {N : A} (α : Path M N) → C N α
path-induction _ b id  =  b
```

This function is defined by *pattern-matching*: we give cases for all possible constructors for $\alpha$, which in this case is just id. In each case, the type of the right-hand side is specialized to that constructor; in this case, C N $\alpha$ is specialized to C M id. It is crucial that path-induction only applies to a family C that is parametrized by a variable x standing for the second endpoint—otherwise, it would imply that all loops of type Path M M are generated by id, which would be incompatible with the homotopy-theoretic interpretation.[3] Topologically, the intuition is that a path with a *free endpoint* can be retracted back to the other, but a path with two fixed endpoints cannot.

[1]We assume basic familiarity with Agda; see Norell [15] for an introduction. We will comment on some of the more idiosyncratic features, such as curly-braces, which are Agda notation for an implicit parameter: we write Path M N when A can be inferred, or Path{A} M N to explicitly notate it.

[2]Agda is a predicative theory, with explicit universe polymorphism, and ordinarily one would define Path in a universe-polymorphic manner. To avoid cluttering the presentation, we use Agda's (inconsistent) `--type-in-type` flag to suppress universe levels, but the development can be done with universe polymorphism instead.

[3]That all loops are id is normally true in Agda, but the `--without-K` option removes this principle.

Paths have a groupoid structure, with inverses and composition defined as follows:

```
! : {A : Type} {M N : A} → Path M N → Path N M
! id  =  id
_∘_ : {A : Type} {M N P : A} → Path N P → Path M N → Path M P
β ∘ id  =  β
```

The groupoid laws hold only up to homotopy; in type theory, this means they are not definitional equalities, but are witnessed by propositional equalities/paths between paths. For example, we can prove associativity, unit, and inverse laws for ∘ ! and id (omitting ∘-unit-r and !-inv-r, which are symmetric):

```
∘-unit-l : {A : Type} {M N : A} (α : Path M N)
        → Path (id ∘ α) α
∘-unit-l id  =  id
∘-assoc : {A : Type} {M N P Q : A}
        (γ : Path P Q) (β : Path N P) (α : Path M N)
        → Path (γ ∘ (β ∘ α)) ((γ ∘ β) ∘ α)
∘-assoc id id id  =  id
!-inv-l  : {A : Type} {M N : A} (α : Path M N)
        → Path (! α ∘ α) id
!-inv-l id  =  id
```

### B. Dependent Types as Fibrations

Just as types act like groupoids, type families act like indexed families of groupoids. In particular, they vary functorially with paths in the indexing type:

```
transport : {B : Type} (E : B → Type)
        {b1 b2 : B} → Path b1 b2 → (E b1 → E b2)
transport C id  =  λ x → x
```

Logically, transport is a "coercion" by propositional equality. transport is functorial up to homotopy; e.g. there is a path between transport C $(\beta \circ \alpha)$ and transport C $\beta$ o transport C $\alpha$, where o is function composition.

While the category-theoretic viewpoint on a type family E : B → Type is as a functor from B to types, the topological viewpoint is as a *fibration*. Given two spaces $\tilde{E}$ and $B$, a fibration over $B$ is a continuous map $p : \tilde{E} \to B$ such that for any $e \in \tilde{E}$, any path in $B$ starting at $p(e)$ has a lifting to a path in $\tilde{E}$ starting at $e$ (and these liftings are continuous/functorial). $\tilde{E}$ is called the *total space*, while $B$ is called the *base space*.

To make the connection with type theory, it is helpful to consider a slightly different characterization of fibrations: Given a point $b$ in the base space $B$, the *fiber over $b$* is the space of points in the total space $\tilde{E}$ that $p$ maps to $b$ (i.e. the inverse image $p^{-1}(b)$). Any path $\beta$ from $b_1$ to $b_2$ in $B$ induces an operation from the fiber over $b_1$ to the fiber over $b_2$. Namely, given $e \in p^{-1}(b_1)$, we lift $\beta$ starting at $e$ and consider the other endpoint of the resulting path; this is well-defined up to homotopy. These operations respect path composition, inversion, and so on, in a homotopical way, yielding a "functor" $E$ sending each $b \in B$ to its fiber $p^{-1}(b)$. Conversely, from any such functor $E$ we can assemble a total space $\tilde{E}$ and a fibration $p : \tilde{E} \to B$ with the specified fibers and path-lifting functions (for groupoids, this is called the "Grothendieck construction").

In type theory, the functor description $E$ corresponds directly to a dependent type E : B → Type, where the type E b represents the fiber over b, and transport E $\alpha$ represents the operation induced by path lifting. Given such an E, the fibration $p : \tilde{E} \to B$ is modeled by the Σ-type $\Sigma b{:}B.E(b)$ and its first projection fst : $\Sigma b{:}B.E(b) \to B$. Thus total spaces are Σ-types, which we will sometimes write as Σ E.

transport "computes" (up to paths) for each specific dependent type, expressing the definition of path lifting for the corresponding fibration; we need the following two rules:

    transport-Path-right : {A : Type} {M N P : A}
      ($\alpha$' : Path N P) ($\alpha$ : Path M N)
      → Path (transport ($\lambda$ x → Path M x) $\alpha$' $\alpha$) ($\alpha$' ∘ $\alpha$)
    transport-→ : {Γ : Type} (A B : Γ → Type) {$\theta 1$ $\theta 2$ : Γ}
      ($\delta$ : $\theta 1$ ≃ $\theta 2$) (f : A $\theta 1$ → B $\theta 1$)
      → Path (transport ($\lambda$ $\gamma$ → (A $\gamma$) → B $\gamma$) $\delta$ f)
        (transport B $\delta$ o f o (transport A (! $\delta$)))

The former says that transporting with the family Path M - is post-composition of paths; the latter that transporting at A → B is given by pre-composing with transport at A (on the inverse) and post-composing with transport at B. Both are proved by matching the input paths as id and returning id.

## C. Functions are Functorial

Simply-typed functions correspond to functors between homotopy types, and have an action at all levels: a function f : A → B has an action on points, paths, paths between paths, etc. The action on points is ordinary function application (f a : B when a : A), while the action on paths is given by

    ap : {A B : Type} {M N : A}
      (f : A → B) → Path {A} M N → Path {B} (f M) (f N)
    ap f id = id

ap acts functorially on identities and composition of paths, and also on identities and composition of functions.

Dependently typed functions f : (b : B) → E (b) correspond to sections of the fibration E. In type-theoretic terms, such a section is a simply-typed function f' : B → Σ E such that fst o f' = ($\lambda$ x → x) *definitionally*—the first component of the result of f' must be its argument. In one direction, given f as above, we can define f' : B → Σ E by $\lambda$ b → (b, f b).

This correspondence helps us arrive at the the dependently typed analogue of ap; the following discussion is a bit difficult, but it is very important for understanding the induction principle for the circle in Section IV. When f has a dependent function type, and $\beta$ : Path {B} $b_1$ $b_2$, then f $b_1$ and f $b_2$ should still be related by a path, but f $b_1$ : E $b_1$ and f $b_2$ : E $b_2$ live in different fibers, i.e. they have different types. Using the above correspondence between f and f', we have ap f' $\beta$ : Path {Σ E} ($b_1$, f $b_1$) ($b_2$, f $b_2$)—there is a path *in the total space* between ($b_1$, f $b_1$) and ($b_2$, f $b_2$). However, we know a little bit more: because fst o f' = ($\lambda$ x → x), it follows from functoriality of ap that ap fst (ap f' $\beta$) equals $\beta$. Topologically, this means that ap f' $\beta$ *projects down to* $\beta$, or *sits above* $\beta$:

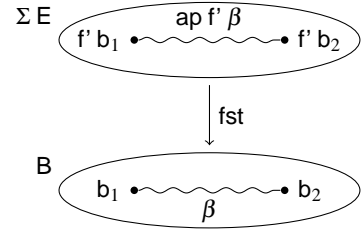

Thus, to describe the result of applying f to $\beta$, it would be intuitively plausible to ask for a path $\tilde{\eta}$ in Σ B such that ap fst $\tilde{\eta}$ *is* $\beta$. However, *is* would need to mean definitional equality here, which cannot be expressed as a type. Instead, we can represent a path between $e_1$ : E $b_1$ and $e_2$ : E $b_2$ sitting above $\beta$ : Path $b_1$ $b_2$ by an $\eta$ : Path {E $b_2$} (transport E $\beta$ $e_1$) $e_2$. That is, we use transport to move $e_1$ into the fiber over $b_2$ and then give a path in E $b_2$. This representation is correct because ($b_1$, $e_1$) is always connected to ($b_2$, transport E $\beta$ $e_1$) by a path over $\beta$ (this follows by path induction on $\beta$), and any path $\tilde{\eta}$ from ($b_1$, $e_1$) to ($b_2$, $e_2$) in Σ E over $\beta$ factors as this path followed by our path $\eta$ in the fiber over $b_2$.

All this motivates the following the dependently typed analogue of ap: applying f to $\beta$ must determine a path between f $b_1$ and f $b_2$ that sits above $\beta$, which is represented by the following type:

    apd : {B : Type} {E : B → Type} {$b_1$ $b_2$ : B}
      (f : (x : B) → E x) ($\beta$ : Path $b_1$ $b_2$)
      → Path (transport E $\beta$ (f $b_1$)) (f $b_2$)
    apd f id = id

## D. Paths Between Functions

Another kind of application is applying a path between functions to an argument, to get a path between results:

    ap≃ : ∀ {A} {B : A → Type} {f g : (x : A) → B x}
      → Path f g → {x : A} → Path (f x) (g x)
    ap≃ $\alpha$ {x} = ap ($\lambda$ f → f x) $\alpha$

    $\lambda$≃ : ∀ {A} {B : A → Type} {f g : (x : A) → B x}
      → ((x : A) → Path (f x) (g x))
      → Path f g

*Function extensionality* $\lambda$≃ is the converse—functions are equal if they are pointwise equal, or a path between functions is a homotopy between them. Agda does not provide this, so we postulate it.[4] We should also give $\beta/\eta$-like equations relating $\lambda$≃ and ap≃, but we do not need them in this paper.

## E. Paths in the Universe

Voevodsky's univalence axiom says roughly that *isomorphic types are equal*, where "equal" means Path, and "isomorphic" means a *homotopy equivalence*,[5] or a pair of mutually inverse functions:

---

[4]It is not strictly necessary to postulate it separately, because it follows from Voevodsky's univalence axiom.

[5]Homotopy equivalences have a slightly undesirable property: there can be multiple different g, $\alpha$, and $\beta$ showing that a given f is a homotopy equivalence. It is common to include an additional coherence cell that fixes this problem. However, any homotopy equivalence can be *improved* by constructing this coherence cell (at the cost of changing $\alpha$ or $\beta$), so we can suppress this detail.

```
record HEquiv (A B : Type) : Type where
  constructor hequiv
  field
    f : A → B
    g : B → A
    α : (x : A) → Path (g (f x)) x
    β : (y : B) → Path (f (g y)) y
```

Rather than stating the univalence axiom in full generality, we specify only the consequences we need. First, a homotopy equivalence determines a path between types:

```
univalence : {A B : Type} → HEquiv A B → Path A B
```

We need two facts about this axiom:

```
transport-univ : {A B : Type} (e : HEquiv A B)
  → Path (transport (λ (A : Type) → A) (univalence e))
         (HEquiv.f e)
!-univalence : {A B : Type} (e : HEquiv A B)
  → Path (! (univalence e))
         (univalence (!-equiv e))
```

The first says that transporting with the identity type family $\lambda\, A \to A$ on an application of univalence applies the forward direction of the equivalence. The second says that the inverse of an application of univalence is the inverse equivalence of e, where !-equiv (hequiv f g $\alpha$ $\beta$) is hequiv g f $\beta$ $\alpha$.

*F. Integers*

We represent integers as follows:

```
data Positive : Type where
  One : Positive
  S : (n : Positive) → Positive
data Int : Type where
  Pos : (n : Positive) → Int
  Zero : Int
  Neg : (n : Positive) → Int
```

It is straightforward to define the successor, predecessor, and addition functions by case-analysis/induction, and to prove that succ and pred are mutually inverse, so we have

```
succ : Int → Int
pred : Int → Int
_+_ : Int → Int → Int
succ-pred : (n : Int) → Path (succ (pred n)) n
pred-succ : (n : Int) → Path (pred (succ n)) n
```

Therefore, we have a homotopy equivalence between Int and itself given by adding and subtracting 1:

```
succEquiv : HEquiv Int Int
succEquiv = hequiv succ pred pred-succ succ-pred
```

### III. INJECTIVITY AND DISJOINTNESS FOR COPRODUCTS

Consider the type of coproducts (binary sums):

```
data _+_ (A B : Type) : Type where
  Inl : A → A + B
  Inr : B → A + B
```

One may expect that constructors are disjoint (Inl a is never equal to Inr b) and injective (Inl a equals Inl a' only if a equals a'). However, a well-known puzzling fact is that this is not provable in pure Martin-Löf type theory, but requires a universe or large eliminations. In this section, we use injectivity and disjointness to introduce the methodology we will use to calculate $\pi_1(S^1)$. The similarities between this proof and our calculation of $\pi_1(S^1)$ offers insight into this puzzling fact, as we explain in Section VI.

Fix A, B, and a:A. Injectivity and disjointness of Inl can be phrased as follows (where Void is the empty type):
- *Injectivity:* If Path (Inl a) (Inl a') then Path a a'.
- *Disjointness:* If Path (Inl a) (Inr b) then Void.

Thus, we can regard the injectivity-and-disjointness problem as the question of characterizing the types Path (Inl a) e for all e. One way to do this is to define a family of types describing the desired characterization, which (for reasons to be explained in Section V) we call Cover:

```
Cover : A + B → Type
Cover (Inl a') = Path a a'
Cover (Inr _) = Void
```

Cover defines a family of types by case analysis, and therefore depends on having a universe Type of types (this is the step that is not possible in pure Martin-Löf type theory). We say that Cover e classifies *codes* for a path in A + B from the fixed *base point* Inl a to the point e.

Suppose that we can encode every path as a code:

```
encode : {e : A + B} → Path (Inl a) e → Cover e
```

Then injectivity and disjointness follow immediately, by the expanding the definition of Cover:

```
inj : {a' : A} → Path (Inl a) (Inl a') → Path a a'
inj {a'} = encode {Inl a'}
dis : {b : B} → Path (Inl a) (Inr b) → Void
dis {b} = encode {Inr b}
```

However, it is easy to define encode, just by transporting along Cover:

```
encode α = transport Cover α id
```

To encode $\alpha$ : Path (Inl a) e, we transport along $\alpha$ with the type family Cover, which reduces the goal of constructing an element of Cover e to that of Cover (Inl a). But Cover (Inl a) is just Path a a, so we can choose id to complete the proof.

Ordinarily, one stops here, having proved injectivity and disjointness, which makes sense when equality is thought of as a mere proposition, without meaningful computational content. However, in homotopy type theory, where paths have real content, it is important to know not only that injectivity and disjointness exist, but that they are equivalences. For example, one might like to know that the paths in the coproduct between Inls are equivalent to the paths in A (i.e. that Path {A + B} (Inl a) (Inl a') is equivalent to Path a a'), so that one may reason about paths in A through their injection into the coproduct. To this end, we will show that encode is an equivalence:

enceqv : {e : A + B} → HEquiv (Path (Inl a) e) (Cover e)
enceqv = hequiv encode decode
              decode-encode encode-decode

by defining decode and proofs decode-encode and encode-decode. decode is defined as follows:

decode : {e : A + B} → Cover e → Path (Inl a) e
decode {Inl a'} α = ap Inl α
decode {Inr _} ()

When e is Inl a', we are given α of type Cover (Inl a'), or Path a a'. Thus, we get a Path (Inl a) (Inl a') by applying Inl to α. When e is Inr, α has type Cover (Inr -), or Void, so the case is vacuously true. We notate this using an *absurd pattern* ().

Next, we show that these two functions are mutually inverse.

### A. Encoding after Decoding

First, we show that starting from a code (an element of the cover), decoding it as a path, and then re-encoding it gives back the original code. We write the proof using a chain of equations, where the notation x ≃⟨ a ⟩ y ≃⟨ b ⟩ z ∎ means there is a path a from x to y and then b from y to z.

encode-decode : {e : A + B} (c : Cover e)
    → encode {e} (decode {e} c) ≃ c
encode-decode {Inl a'} α =
    encode (decode α)          -- (1)
      ≃⟨ id ⟩
    transport Cover (ap Inl α) id    -- (2)
      ≃⟨ ap≃ (! (transport-ap-assoc' Cover Inl α)) ⟩
    transport (Cover o Inl) α id    -- (3)
      ≃⟨ id ⟩
    transport (λ a' → Path a a') α id    -- (4)
      ≃⟨ transport-Path-right α id ⟩
    α ∘ id                      -- (5)
      ≃⟨ id ⟩
    α ∎                         -- (6)
encode-decode {Inr _} ()

The proof begins by casing on e. In the Inl a' case, α is a path from a to a'. Between line 1 and line 2, we expand the definitions of encode and decode, where for decode we know the case for Inl is selected. Between line 2 and line 3, we reassociate transport and ap: in general, transport (C o f) α is the same as transport C (ap f α). Between lines 3 and 4, we reduce the definition of transport on Inl. Between lines (4) and (5), we apply the fact that transporting at Path a - is post-composition. Between lines 5 and 6, we apply one of the unit laws for composition, which gives the result.

The Inr case is vacuously true, because in this case we have an element of Cover (Inr -), which is the empty type.

### B. Decoding after Encoding

The other direction is

decode-encode : {e : A + B} (α : Path (Inl a) e)
    → Path (decode {e} (encode {e} α)) α

Expanding the defnition of encode, we need a path from (decode (transport Cover α id)) to α, where α is an arbitrary path from Inl (a) to e. The key idea is to apply *path induction*: To prove the goal for an arbitrary e : A + B and α : Path (Inl a) e, it suffices to consider the case where e is Inl a and α is id. In this case, we have to show that decode (encode id) is id, which is easy: it holds definitionally, because both transport and ap compute to id on id, so id proves the result. This argument is formalized as follows:

decode-encode {e} α =
    path-induction
      (λ e' α' → Path (decode {e'} (encode {e'} α')) α')
      id α

### C. Summary

To review, the structure of this proof is: (1) Fix a base point, and define the "cover", which is a type of codes for paths from the base point. (2) Define encode by transporting along the cover, with an appropriate base case. (3) Define decode by case-analysis/induction. (4) Prove that encoding after decoding is the identity, using the induction principle for codes. (5) Prove that decoding after encoding is the identity, using path induction. We will follow this same template for the circle.

## IV. THE CIRCLE

In this section, we introduce the representation of the circle in homotopy type theory. The natural numbers is an inductive type, with *generators* zero and successor. One of the new ingredients in homotopy type theory is *higher-dimensional inductive types* (or just *higher inductive types*) [12, 13, 17]: inductive types specified by generators not only for points (terms), but also for paths.

One might draw the circle like this:

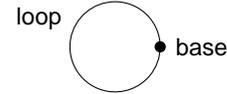

It has a single point, and a single non-identity loop from this base point to itself. This translates to a higher inductive type with two generators:

base : $S^1$
loop : Path {$S^1$} base base

base is like an ordinary constructor for an inductive type, with the same status as zero or successor. loop is similar, except it generates a path on the circle. This generator can be used with the generic groupoid structure to form additional paths, such as ! loop, loop ∘ loop, etc. Some of these paths are "reducible"— for example, loop ∘ ! loop is homotopic to id.

### A. Simple Elimination

That the type of natural numbers is *inductively* generated by zero and successor is expressed by its elimination rule: to map from the natural numbers into a type X, it suffices to specify an element and an endomorphism of X, to be the images of zero and successor. Similarly, the elimination rule for $S^1$ expresses that the circle is inductively generated by base and loop: to map from the circle into any other type, it suffices to find a point and a loop in that type:

```
S¹-recursion : {X : Type}
               (base' : X) (loop' : Path base' base')
               → S¹ → X
```

Elimination rules require accompanying $\beta$-reduction rules, which for an inductive type state that "the elimination rule, applied to a generator, computes to the corresponding branch". In this case, this means that S¹-recursion should compute to base' when applied to base and to loop' when applied to loop:

```
(S¹-recursion base' loop') base = base'
(S¹-recursion base' loop') loop = loop'
```

However, the second equation does not quite make sense, because S¹-recursion base' loop' is a function S¹ → X and loop is a path in S¹. Thus, we need to use ap to apply this function to the path:

```
ap (S¹-recursion base' loop') loop = loop'
```

The left-hand side is a path from (S¹-recursion base' loop') base to itself, which by the first $\beta$-reduction is Path base' base'; so the two sides have the same type.

Thus, the computation rules follow the same general pattern as for ordinary inductive types, except we need to use the notion of application appropriate for the level of the generator.

*B. Dependent Elimination*

To fully characterize an inductive type, we need not just recursion, but an induction principle (dependent elimination). The induction principle for natural numbers says that to prove a property of natural numbers (i.e. inhabit a type family indexed over natural numbers), it suffices to prove that it holds for zero and is preserved by successors. Similarly, the induction rule for S¹ says that to prove a property of points on the circle, it suffices to prove that it holds for base and is "preserved by going around the loop".

```
S¹-induction : (X : S¹ → Type)
               (base' : X base)
               (loop' : Path (transport X loop base') base')
               → (y : S¹) → X y
```

Here, X is not just a type, but a dependent type/fibration over the circle. The natural thing to ask is that base', the image of base, should show that X holds for base, or be a point in the fiber over base. loop' must be a path from base' to itself, but one that *projects down to loop* in the sense described in Section II-C—recall that a path from base' to itself that sits above loop in the fibration X is represented by the type given to loop' above. To see that this is appropriate, note that the result type (y : S¹) → X y represents topologically a *section* of the fibration $\Sigma X \to S^1$, which must take each point or path in S¹ to a point or path lying above it; thus we should expect to need a path lying above loop as input. More syntactically, the point is basically that the second computation rule for S¹-induction must be well-typed:

```
(S¹-induction base' loop') base = base'
apd (S¹-induction base' loop') loop = loop'
```

The dependent elimination rule also characterizes an inductive type up to equivalence. S¹-induction is equivalent to asserting that for all X, the type of functions S¹ → X is naturally equivalent to the type of pairs (base' : X, loop' : Path base' base') (the premises of S¹-recursion). Because the type theory ensures that functions are groupoid homomorphisms, this is exactly the universal property of the free ∞-groupoid generated by one object and one loop on it.

*C. Agda Implementation*

Computer proof assistants such as Agda and Coq include ordinary inductive types, but not yet higher inductive types. One way to implement the latter is to simply postulate the generators, the elimination rule, and the computation rules. With this representation, the computation rules are paths (propositional equalities), rather than definitional equalities.

Though the question of whether these rules should be definitional equalities or paths has not yet been settled, it is certainly more convenient to do proofs if they are definitional equalities. While it is not possible to achieve this in Agda, there is a trick using private data types that allows the base (but not loop) rule to be definitional [8]. Because this simplifies the proofs, we use this implementation, with a postulate

```
βloop/rec : {X : Type}
            (base' : X) (loop' : Path base' base')
            → Path (ap (S¹-recursion base' loop') loop) loop'
```

for the $\beta$-rule for loop (and similarly for S¹-induction).

V. THE FUNDAMENTAL GROUP OF THE CIRCLE

Given a space $X$ with a specified base point $x_0$, the fundamental group $\pi_1(X, x_0)$ (or just $\pi_1(X)$ when $x_0$ is clear from context) is the group of homotopy classes of loops from $x_0$ to itself, with path composition as the group operation. In type theory, this corresponds to the type Path {X} $x_0$ $x_0$, except for one caveat: For $\pi_1(X)$, the group has a *set* of elements, which are paths *quotiented by homotopy*. This means that any two paths that are homotopic are equal, but any non-trivial *structure* of paths between paths has been collapsed by quotienting. On the other hand, the *type* Path $x_0$ $x_0$ may have interesting paths between paths (i.e., the type Path {Path $x_0$ $x_0$} $\alpha$ $\beta$ might not be trivial). Thus, Path $x_0$ $x_0$ corresponds more closely to what is called the *loop space* $\Omega_1(X, x_0)$, the *space* of loops in $X$ based at $x_0$, which also may still have non-trivial structure.

It is possible to construct the set $\pi_1(X)$ from the loop space $\Omega_1(X)$ by an operation called *truncation* (in classical topology, this is just the set of connected components). However, for the example we consider here, this is not necessary: what we will prove is that the loop space of the circle $\Omega_1(S^1)$ is $\mathbb{Z}$. Because the truncation of $\mathbb{Z}$ is $\mathbb{Z}$, applying truncation to both sides shows that $\pi_1(S^1)$ is also $\mathbb{Z}$. Moreover, proving that $\Omega_1(S^1)$ is $\mathbb{Z}$ immediately characterizes all the *higher homotopy groups* of the circle, because the higher homotopy groups are determined by iterating the loop space construction. Thus, if $\Omega_1(S^1)$ is $\mathbb{Z}$, then the higher homotopy groups of $S^1$ must be the same as the higher homotopy groups of $\mathbb{Z}$, and therefore trivial.

In type theoretic terms, this means that our first goal is to prove that the type Path {$S^1$} base base is equivalent to Int. We then check that this equivalence is a group homomorphism, taking path composition to addition.

### A. Classical and Type-theoretic Proofs

Our proof that $\Omega_1(S^1)$ is $\mathbb{Z}$ can be seen as a type-theoretic version of a proof in classical homotopy theory. The classical proof we start from is usually formulated using "universal covering spaces", but we will give an equivalent sketch using *fibrations* which transfers more directly to type theory. Recall the notion of a fibration from Section II-B. For any point $x_0 \in B$, there is a canonical *path fibration* $p : P_{x_0}B \to B$, where the points of $P_{x_0}B$ are paths in $B$ starting at $x_0$, and the map $p$ selects the other endpoint of such a path. The space $P_{x_0}B$ is *contractible*, i.e. homotopy equivalent to a point, since we can "retract" any path to its initial endpoint $x_0$. Moreover, the fiber over $x_0$ is the loop space $\Omega_1(B, x_0)$.

Now consider the "winding" map $w : \mathbb{R} \to S^1$, which looks like a helix projecting down onto the circle:

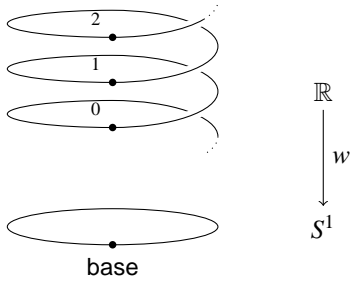

The map $w$ sends each point on the helix to the point on the circle that it is "sitting above"; this map is a fibration, and the fiber over each point is isomorphic to the integers. If we lift the path that goes counterclockwise around the loop on the bottom, we go up one level in the helix, incrementing the integer in the fiber. Similarly, going clockwise around the loop on the bottom corresponds to going down one level in the helix, decrementing this count. This fibration is called the *universal cover* of the circle.

Now a basic fact is that a map $E_1 \to E_2$ of fibrations over $B$ which is a homotopy equivalence between $E_1$ and $E_2$ induces a homotopy equivalence on fibers. Since $\mathbb{R}$ and $P_{\text{base}}S^1$ are both contractible, they are homotopy equivalent, and thus the fibers over base, $\mathbb{Z}$ and $\Omega_1(S^1)$, are isomorphic.

It is possible to formalize this classical proof directly in type theory, by (1) defining the universal cover, (2) proving that a homotopy equivalence between total spaces induces an equivalence on fibers, and (3) proving that the total spaces of both the path fibration and the cover are contractible (this was the first proof of this result in homotopy type theory [18]). However, it turns out to be simpler to explicitly construct the encoding-decoding equivalence, following the template introduced in Section III [9]. Both proofs use the same construction of the cover (step 1 above). Where the classical proof induces an equivalence on fibers from an equivalence between total spaces (step 2), the type-theoretic proof constructs the inverse map explicitly as a map between fibers. Where the classical proof uses contractibility (step 3), the type-theoretic proof uses path induction, circle induction, and integer induction. These are the same tools used to prove contractibility—indeed, path induction *is* contractibility of the path fibration composed with transport—but it is more convenient to use them to prove inverses directly. This proof is a good example of how combining insights from homotopy theory and type theory can simplify proofs and yield deeper insight.

### B. The Universal Cover of the Circle

Recall that fibrations are represented by dependent types (Section II-B). The path fibration $P_{\text{base}}S^1 \to S^1$ is easy to represent: it is the dependent type sending x : $S^1$ to Path base x. The universal cover of the circle (the helix) requires a bit more thought: since our "circle" is not a topological one, we don't have a "real line" that we can wrap around it. However, our inductive definition of the circle gives us a different way to define dependent types over it: by circle-recursion.

```
Cover : S¹ → Type
Cover x = S¹-recursion Int (univalence succEquiv) x
```

To define a function by circle recursion, we need to find a point and a loop in the target. In this case, the target is Type, and the point we choose is Int, corresponding to our expectation that the fiber of the universal cover should be the integers. The loop we choose is the successor/predecessor isomorphism on Int, succEquiv (Section II), which by univalence determines a path from Int to Int. Univalence is necessary for this part of the proof, because we need a *non-trivial* loop from Int to Int.

It is immediate from this definition that Cover base is Int, and we can verify that choosing succEquiv as the image of loop gives the desired path-lifting action: transporting one way along the cover is successor (going up one level in the helix), and the other way is predecessor (going down one level):

```
transport-Cover-loop : Path (transport Cover loop) succ
transport-Cover-loop =
    transport Cover loop
    ≃⟨ transport-ap-assoc Cover loop ⟩
    transport (λ x → x) (ap Cover loop)
    ≃⟨ ap (transport (λ x → x))
         (βloop/rec Int (univalence succEquiv)) ⟩
    transport (λ x → x) (univalence succEquiv)
    ≃⟨ transport-univ _ ⟩
    succ ∎

transport-Cover-!loop : Path (transport Cover (! loop)) pred
```

For transport-Cover-loop, we re-associate, which creates a $\beta$-redex ap Cover loop for $S^1$-recursion. After reducing this, we have a term of of the form transport (λ x → x) (univalence e), which is a $\beta$-redex for univalence, and selects the "forward" direction of the equivalence. We omit the proof for transport-Cover-!loop, which is similar except for additional reasoning about inverses, using !-univalence from Section II.

The cover can be seen as a type of codes for paths on the circle. The next step is to define "encoding" and "decoding"

functions and prove that they are an equivalence between Path base x and Cover x for all x : $S^1$. This is a generalization of the original statement we intended to prove, which was that Path base base is equivalent to Cover base (the latter being, by definition, Int).

*C. Encoding*

As in Section III, encode is defined by transporting in the cover. The starting point needs to be an element of Cover base, which is Int. In this case, a good choice is Zero[6]:

 encode : {x : $S^1$} → Path base x → Cover x
 encode $\alpha$ = transport Cover $\alpha$ Zero

The instance encode' = encode {base} has type Path base base → Int, as we originally intended.

The interesting thing about this function is that it computes a concrete number from a loop on the circle, when this loop is represented using the abstract groupoidal framework of homotopy type theory. To gain an intuition for how it does this, observe that by the above lemmas, transport Cover loop is succ and transport Cover (! loop) is pred. Further, transport is functorial (Section II), so transport Cover (loop ∘ loop) is (transport Cover loop) o (transport Cover loop), etc. Thus, when $\alpha$ is a composition like

 loop ∘ ! loop ∘ loop ∘ ...

transport Cover $\alpha$ will compute a composition of functions like

 succ o pred o succ o ...

Applying this composition of functions to Zero will compute the *winding number* of the path—how many times it goes around the circle, with orientation marked by whether it is positive or negative, after inverses have been canceled.

Thus, the computational content of encode follows from the $\beta$-like rules for higher-inductive types and univalence, and the action of transport on compositions and inverses. This "computation" happens only up to paths in current homotopy type theory, so we cannot actually run this program in Agda, but an alternate formulation might take these equations as definitional equalities [10].

*D. Decoding*

The first step in decoding is that, given an integer *n*, we compute the *n*-fold composition loop$^n$:

 loop^ : Int → Path base base
 loop^ Zero     = id
 loop^ (Pos One) = loop
 loop^ (Pos (S n)) = loop ∘ loop^ (Pos n)
 loop^ (Neg One) = ! loop
 loop^ (Neg (S n)) = ! loop ∘ loop^ (Neg n)

At this point, if we were working naively, rather than with Section III or the classical proof in mind, we might think that this is enough. That is, since what we want overall is an equivalence between Path base base and Int, we might expect

---

[6]We could choose another number as the base case besides Zero, but then we would need to subtract it off when decoding.

to be able to prove that encode' : Path base base → Int and loop^ : Int → Path base base give an equivalence. The problem comes in trying to prove the "decode after encode" direction:

 decode-encode : {$\alpha$ : Path base base}
      → Path (loop^ (encode' $\alpha$)) $\alpha$

In Section III, we proved this step using path induction to reduce $\alpha$ to the identity, which depends crucially on $\alpha$ having one endpoint free—recall that path induction does not apply to loops like a Path base base with both endpoints fixed! The way to solve this problem is to state decode-encode generally for all x:$S^1$ and $\alpha$ : Path base x:

 decode-encode : {x:$S^1$} {$\alpha$ : Path base x}
      → Path (loop^ (encode {x} $\alpha$)) $\alpha$

However, this does not type check as is, because loop^ works only for Path base base, whereas here we need Path base x. This gives a direct way to see the necessity of extending loop^ to a function with a more general type:

 decode : {x : $S^1$} → Cover x → Path base x

Of course, the template of Section III and the proof from classical homotopy theory also lead us to expect to need such a generalization.

Here is the definition of decode:

 decode : {x : $S^1$} → Cover x → Path base x
 decode {x} =
 $S^1$-induction
   ($\lambda$ x' → Cover x' → Path base x')
   loop^
   ( transport ($\lambda$ x' → Cover x' → Path base x') loop loop^
   $\simeq$ transport ($\lambda$ x' → Path base x') loop
       o loop^
       o transport Cover (! loop)
   $\simeq$ ($\lambda$ p → loop ∘ p) o loop^ o transport Cover (! loop)
   $\simeq$ ($\lambda$ p → loop ∘ p) o loop^ o pred
   $\simeq$ ($\lambda$ n → loop ∘ (loop^ (pred n)))
   $\simeq$ ($\lambda$ n → loop^ n)
   ∎)
 x

decode's first argument is an arbitrary point on the circle. Thus, we proceed by circle induction, which requires (1) a function Cover base → Path base base, which is just loop^, and (2) a path showing that this function is preserved by going around the loop. Formally, this means a path from transport ($\lambda$ x' → Cover x' → Path base x') loop loop^ to loop^.

Above, we have shown the steps of reasoning required to give such a path, eliding the proof terms, which we now discuss informally. The path is constructed by composing five steps of reasoning, between the six lines above, starting from transport ($\lambda$ x' → Cover x' → Path base x') loop loop^. From line 1 to line 2, we apply the definition of transport when the outer connective of the type family is →, using the lemma transport-→ from Section II. This reduces the transport to pre- and post-composition with transport at the domain and range types. From line 2 to line 3, we apply the definition of transport when the type family is Path base - (called transport-Path-right above). From line 3 to line 4, we apply transport-Cover-!loop. From line 4 to line 5, we simply reduce the function composition. The final step is the only significant one: it follows from associativity and inverses of ∘, together with a lemma loop^-preserves-pred which gives a Path (loop^ (pred n)) (! loop ∘ loop^ n) for all n. This lemma is proved by a simple case analysis, again using associativity/unit/inverse laws.

### E. Encoding after Decoding

Computing encode ∘ loop^ is comparatively straightforward.

```
encode-loop^ : (n : Int) → Path (encode (loop^ n)) n
encode-loop^ Zero = id
encode-loop^ (Pos One) = ap≃ transport-Cover-loop
encode-loop^ (Pos (S n)) =
  encode (loop^ (Pos (S n)))
    ≃⟨ id ⟩
  transport Cover (loop ∘ loop^ (Pos n)) Zero
    ≃⟨ ap≃ (transport-∘ Cover loop (loop^ (Pos n))) ⟩
  transport Cover loop
          (transport Cover (loop^ (Pos n)) Zero)
    ≃⟨ ap≃ transport-Cover-loop ⟩
  succ (transport Cover (loop^ (Pos n)) Zero)
    ≃⟨ id ⟩
  succ (encode (loop^ (Pos n)))
    ≃⟨ ap succ (encode-loop^ (Pos n)) ⟩
  succ (Pos n) ∎
```

The proof is a simple induction on Int, using functoriality of transport and the transport-Cover-loop lemmas. We omit the cases for Neg, which are analogous.

To prove the equivalence of Path base base and Int, this is sufficient. If we additionally want an equivalence between Path base x and Cover x for general x, then we need to show that encode-loop^ extends to

```
encode-decode : {x : S¹} → (c : Cover x)
                → Path (encode (decode {x} c)) c
encode-decode {x} = S¹-induction
  (λ (x : S¹) → (c : Cover x)
                → Path (encode {x} (decode {x} c)) c)
  encode-loop^ proof x
```

This can be proved using circle induction, with encode-loop^ as the image of base. The proof that encode-loop^ is compatible with the loop requires a path between two paths in Int. The easiest way to define this is to observe that all paths between paths in Int are equal (that is, Int is an *hset* or satisfies uniqueness of identity proofs), which can be proved by showing that it has decidable equality and then applying Hedberg's theorem [5].

### F. Decoding after Encoding

As in Section III, the proof for decoding after encoding is a single path induction: Suppose α is id : Path base base. Then encode {base} id = Zero, and decode {base} Zero = loop^ Zero = id, so we need a Path id id—which can be id.

```
decode-encode : {x : S¹} (α : Path base x)
                → Path (decode (encode α)) α
decode-encode {x} α =
path-induction
  (λ (x' : S¹) (α' : Path base x')
       → Path (decode (encode α')) α')
  id α
```

decode-encode can be seen as an η-rule/induction principle for paths on the circle, which states that every loop on the circle is of the form loop^ n for some n:

```
all-loops : (α : Path base base) → Path α (loop^ (encode α))
all-loops α = ! (decode-encode α)
```

Consequently, to prove a statement for every Path base base, it suffices to prove the statement for every path loop^ n, which can be done using integer induction on n. Recall that the proof of decode-encode depends crucially on the fact that decode is defined for a path with a free endpoint. Thus, the essential parts of this proof are the path induction used here, and the circle induction used to define decode from loop^. It is crucial to the methodology for working in homotopy type theory that this combination of circle and path induction suffices to *prove* this induction principle for loops on the circle, as we discuss further in Section VI.

### G. Summary

These lemmas give a homotopy equivalence between Path base base and Int, establishing that the loop space of the circle is equivalent to Int:

```
Ω₁ [S¹] -is-Int : HEquiv (Path base base) Int
Ω₁ [S¹] -is-Int =
    hequiv encode decode decode-encode encode-loop^
```

To identify the fundamental group of $S^1$ with Int as a *group*, we also must check that this equivalence is a group homomorphism. A bijection between carriers is a group homomorphism if one of the functions preserves composition (it then necessarily preserves inverses and the unit because these are unique). Thus, it suffices to show that

```
preserves-composition : (n m : Int)
    → Path (loop^ (n + m)) (loop^ n ∘ loop^ m)
```

The proof is an easy induction, using associativity and unit of ∘, the lemma loop^-preserves-pred defined above, and an analogous lemma that loop^ (succ n) is loop ∘ loop^ n.

Categorically, we can understand the proof we have just given as follows: The higher-inductive type $S^1$ is a description of the free ∞-groupoid with one morphism, represented *abstractly* using the groupoidal framework of type theory. The proof we have given shows that type theory is sufficiently powerful to relate this abstract description to a *concrete* description of the free group on one generator, as the inductive type Int equipped with the function +.

## VI. CONCLUSION

In this paper, we have described a technique for characterizing the path spaces of inductive types in type theory, and applied it to two examples. For coproducts, we obtain injectivity and disjointness of constructors. For the circle, we compute its fundamental group, a basic theorem of algebraic topology. The proof for the circle illustrates the use of homotopy type theory as a logic of homotopy theory: using higher inductive types and the ambient groupoidal framework of the type theory, we can represent homotopy types and prove interesting mathematical properties of them. Our technique extends to other types: Kuen-Bang Hou (Favonia), Chris Kapulkin, Carlo Angiuli, and the first author have used the same methodology to prove that the fundamental group of a bouquet of $n$ circles ($n$ circles around a single point) is the free group on $n$ generators.

Seeing injectivity-and-disjointness in this context provides a topological explanation for the use of a universe to prove them: Injectivity and disjointness characterize a path space. Topological proofs characterizing a path space typically consider an entire path fibration at once (like Path (Inl a) −), rather than a path with both endpoints fixed (like Path (Inl a) (Inl a′)), and show that the entire path fibration is equivalent to an alternate fibration (our "codes"). The codes fibration (like the universal cover of the circle) is represented in type theory using induction, which requires a universe or large elimination.

Moreover, the fact that a universe is *necessary* has analogues in higher dimensions: it is the first rung on a ladder of categorical nondegeneracy. Without a universe, the category of types could be a poset, in which case disjointness at least would fail. Without a *univalent* universe, the ∞-category of types could be a 1-category, in which case the computation of the fundamental group of the circle would fail. In general, path spaces of inductive types are only "correct" when the category of types is sufficiently rich to support them.

In this paper, we have taken an approach to inductive types where the characterization of the path space is a *theorem*, not part of the *definition*. One might wonder whether we could take the opposite approach: For coproducts, we might include injectivity and disjointness of Inl and Inr in the definition; for the circle, we might include an elimination rule for paths Path {S¹} x y expressing that they are freely generated by loop. However, there are two problems with this. Conceptually, a (higher) inductive type is *one* freely generated structure, even though it may have more than one kind of generator. As such, it should have only one elimination rule, expressing its universal property. More practically, calculating homotopy groups of a space in algebraic topology can be a significant mathematical theorem. For example, for the two-dimensional sphere, $\pi_1$ is trivial, $\pi_2$ is $\mathbb{Z}$ (like the circle, one level up), but $\pi_3$ is also $\mathbb{Z}$, *even though the description of the sphere does not include any generators at this level*. This is due to something called the Hopf fibration, which arises from the interaction of the lower-dimensional generators with the ∞-groupoid laws. Indeed, there is no general formula known for the homotopy groups of higher-dimensional spheres, so we would not know what characterization to include in the definition, even if we wanted to.

Fortunately, the examples in this paper suggest that the paths in inductive types will always be *determined* by the inductive description and ambient ∞-groupoid laws—so characterizing the path spaces explicitly in the definition would be at best redundant, and at worst inconsistent. Thus, we can *pose* these questions about homotopy groups using higher inductive types, and hope to use homotopy type theory to answer them.

**Acknowledgments** We thank Steve Awodey, Robert Harper, Chris Kapulkin, Peter Lumsdaine, and many other people in Carnegie Mellon's HoTT group and the Institute for Advanced Study special year for helpful discussions about this work.

## References

[1] S. Awodey and M. Warren. Homotopy theoretic models of identity types. *Mathematical Proceedings of the Cambridge Philosophical Society*, 2009.

[2] Coq Development Team. *The Coq Proof Assistant Reference Manual, version 8.2*. INRIA, 2009. Available from http://coq.inria.fr/.

[3] N. Gambino and R. Garner. The identity type weak factorisation system. *Theoretical Computer Science*, 409(3):94–109, 2008.

[4] R. Garner. Two-dimensional models of type theory. *Mathematical. Structures in Computer Science*, 19(4):687–736, 2009.

[5] M. Hedberg. A coherence theorem for Martin-Löf's type theory. *Journal of Functional Programming*, 8(4):413–436, July 1998.

[6] M. Hofmann and T. Streicher. The groupoid interpretation of type theory. In *Twenty-five years of constructive type theory*. Oxford University Press, 1998.

[7] C. Kapulkin, P. L. Lumsdaine, and V. Voevodsky. The simplicial model of univalent foundations. arXiv:1211.2851, 2012.

[8] D. R. Licata. Running circles around (in) your proof assistant; or, quotients that compute. http://homotopytypetheory.org/2011/04/23/running-circles-around-in-your-proof-assistant/, April 2011.

[9] D. R. Licata. A simpler proof that $\pi_1(S^1)$ is $\mathbb{Z}$. http://homotopytypetheory.org/2012/06/07/a-simpler-proof-that-$\pi_1$s¹-is-z/, June 2012.

[10] D. R. Licata and R. Harper. Canonicity for 2-dimensional type theory. In *ACM SIGPLAN-SIGACT Symposium on Principles of Programming Languages*, 2012.

[11] P. L. Lumsdaine. Weak ω-categories from intensional type theory. In *International Conference on Typed Lambda Calculi and Applications*, 2009.

[12] P. L. Lumsdaine. Higher inductive types: a tour of the menagerie. http://homotopytypetheory.org/2011/04/24/higher-inductive-types-a-tour-of-the-menagerie/, April 2011.

[13] P. L. Lumsdaine and M. Shulman. Higher inductive types. In preparation, 2013.

[14] U. Norell. *Towards a practical programming language based on dependent type theory*. PhD thesis, Chalmers University of Technology, 2007.

[15] U. Norell. Dependently typed programming in agda. *Summer school on Advanced Functional Programming*, 2008.

[16] C. Paulin-Mohring. *Extraction de programmes dans le Calcul des Constructions*. PhD thesis, Université Paris 7, 1989.

[17] M. Shulman. Homotopy type theory VI: higher inductive types. http://golem.ph.utexas.edu/category/2011/04/homotopy_type_theory_vi.html, April 2011.

[18] M. Shulman. A formal proof that $\pi_1(S^1) = \mathbb{Z}$. http://homotopytypetheory.org/2011/04/29/a-formal-proof-that-pi1s1-is-z/, April 2011.

[19] B. van den Berg and R. Garner. Types are weak ω-groupoids. *Proceedings of the London Mathematical Society*, 102(2):370–394, 2011.

[20] V. Voevodsky. Univalent foundations of mathematics. Invited talk at WoLLIC 2011 18th Workshop on Logic, Language, Information and Computation, 2011.

[21] M. A. Warren. *Homotopy theoretic aspects of constructive type theory*. PhD thesis, Carnegie Mellon University, 2008.